\documentclass[11pt]{amsart}
\usepackage{amsmath,amssymb,mathrsfs}
\newtheorem{theorem}{Theorem}[section]
\newtheorem{proposition}[theorem]{Proposition}
\newtheorem{corollary}[theorem]{Corollary}
\newtheorem{lemma}[theorem]{Lemma}
\begin{document}

\title[Constant mean curvature surfaces]{Constant mean curvature surfaces in warped product manifolds}
\author{Simon Brendle}
\thanks{The author was supported in part by the National Science Foundation under grant DMS-0905628.}
\begin{abstract}
We consider surfaces with constant mean curvature in certain warped product manifolds. We show that any such surface is umbilic, provided that the warping factor satisfies certain structure conditions. This theorem can be viewed as a generalization of the classical Alexandrov theorem in Euclidean space. In particular, our results apply to the deSitter-Schwarzschild and Reissner-Nordstrom manifolds.
\end{abstract}
\address{Department of Mathematics \\ Stanford University \\ Stanford, CA 94305}
\maketitle

\section{Introduction} 

A classical theorem due to Alexandrov \cite{Alexandrov} asserts that any closed, embedded hypersurface in $\mathbb{R}^n$ with constant mean curvature is a round sphere. Alexandrov's theorem is remarkable in that it holds in all dimensions; requires no assumptions about the topology of the surface; and does not impose any stability assumptions. More generally, it is known that surfaces 
of constant mean curvature in the hemisphere and in hyperbolic space are geodesic spheres (see e.g. \cite{Hijazi-Montiel-Roldan}, \cite{Montiel-Ros}). Montiel \cite{Montiel} has obtained a uniqueness theorem for star-shaped hypersurfaces of constant mean curvature in certain rotationally symmetric manifolds. The argument in Montiel's paper \cite{Montiel} applies to various ambient spaces; the assumption that the surface is star-shaped plays a crucial role in the argument (see also \cite{Barbosa-doCarmo}).

In a different direction, Christodoulou and Yau \cite{Christodoulou-Yau} studied stable surfaces of constant mean curvature in asymptotically flat three-manifolds. Their work was motivated by considerations in general relativity; in particular, they showed that any such surface has nonnegative Hawking mass. Later, Bray \cite{Bray-thesis} studied the isoperimetric problem in the three-dimensional Schwarzschild manifold. A surface is called isoperimetric if it has minimal area among all surfaces that enclose the same volume. Using an ingenious comparison argument, Bray \cite{Bray-thesis} was able to show that any isoperimetric surface in the Schwarzschild manifold must be a coordinate sphere. We note that Bray's method can be extended to more general ambient manifolds; see \cite{Bray-Morgan} and \cite{Corvino-Gerek-Greenberg-Krummel} for details. 

In 1997, Huisken and Yau \cite{Huisken-Yau} proved that for any asymptotically flat three-manifold $M$ with positive ADM mass, there exists a compact set $K$ so that the complement $M \setminus K$ can be foliated by stable surfaces of constant mean curvature. Moreover, Huisken and Yau proved a uniqueness result for such foliations near infinity under mild additional assumptions. Qing and Tian \cite{Qing-Tian} have obtained a uniqueness result for stable surfaces of constant mean curvature that avoid some large compact set $K$. It was shown by Eichmair and Metzger \cite{Eichmair-Metzger2} that a stable constant mean curvature surface $\Sigma$ must avoid a given compact set $K$, provided that the area of $\Sigma$ is sufficiently large and the ambient manifold $M$ has positive scalar curvature. Eichmair and Metzger also proved that the constant mean curvature spheres constructed by Huisken and Yau are, in fact, isoperimetric surfaces. This confirmed a conjecture of Bray; see \cite{Eichmair-Metzger1} and \cite{Eichmair-Metzger3} for details. Finally, we note that Rigger \cite{Rigger} and Neves and Tian \cite{Neves-Tian} have constructed foliations by surfaces of constant mean curvature in asymptotically hyperbolic manifolds.

In this paper, we prove an analogue of Alexandrov's theorem for a class of warped product manifolds. Let us fix an integer $n \geq 3$. Throughout this paper, we assume that $N$ is a compact Riemannian manifold of dimension $n-1$ such that 
\[\text{\rm Ric}_N \geq (n-2)\rho \, g_N\] 
for some constant $\rho$. Moreover, we consider a smooth positive function $h: [0,\bar{r}) \to \mathbb{R}$ which satisfies the following conditions: 
\begin{itemize} 
\item[(H1)] $h'(0) = 0$ and $h''(0) > 0$.
\item[(H2)] $h'(r) > 0$ for all $r \in (0,\bar{r})$.
\item[(H3)] The function 
\[2 \, \frac{h''(r)}{h(r)} - (n-2) \, \frac{\rho - h'(r)^2}{h(r)^2}\] 
is non-decreasing for $r \in (0,\bar{r})$.
\item[(H4)] We have 
\[\frac{h''(r)}{h(r)} + \frac{\rho-h'(r)^2}{h(r)^2} > 0\] 
for all $r \in (0,\bar{r})$.
\end{itemize}
We now consider the manifold $M = N \times [0,\bar{r})$ equipped with the Riemannian metric 
\begin{equation} 
\label{metric}
g = dr \otimes dr + h(r)^2 \, g_N. 
\end{equation}
The following is the main result of this paper:

\begin{theorem} 
\label{alexandrov}
Suppose that $(M,g)$ is a warped product manifold satisfying conditions (H1)--(H3). Moreover, let $\Sigma$ be a closed, embedded, orientable hypersurface in $(M,g)$ with constant mean curvature. Then $\Sigma$ is umbilic. If, in addition, the condition (H4) holds, then $\Sigma$ is a slice $N \times \{r\}$ for some $r \in (0,\bar{r})$.
\end{theorem} 

It is interesting to consider the special case when $\text{\rm Ric}_N = (n-2)\rho \, g_N$. In this case, the Ricci and scalar curvature of $g$ are given by 
\begin{align} 
\label{ricci.tensor}
\text{\rm Ric} 
&= -\Big ( \frac{h''(r)}{h(r)} - (n-2) \, \frac{\rho-h'(r)^2}{h(r)^2} \Big ) \, g \notag \\ 
&- (n-2) \, \Big ( \frac{h''(r)}{h(r)} + \frac{\rho-h'(r)^2}{h(r)^2} \Big ) \, dr \otimes dr 
\end{align} 
and 
\begin{equation} 
R = -(n-1) \, \Big ( 2 \, \frac{h''(r)}{h(r)} - (n-2) \, \frac{\rho - h'(r)^2}{h(r)^2} \Big ). 
\end{equation}
Hence, in this case, the condition (H3) is equivalent to saying that the scalar curvature of $g$ is non-increasing in $r$. Moreover, the condition (H4) says that the Ricci curvature is smallest in the radial direction.

In particular, condition (H3) is satisfied if $N$ is the standard sphere and $(M,g)$ has constant scalar curvature. Besides the standard spaces of constant sectional curvature, the most basic examples of rotationally symmetric manifolds with constant scalar curvature are the deSitter-Schwarzschild manifolds. We briefly recall their definition. Let us fix real numbers $m$ and $\kappa$. We assume that $m$ is positive. Moreover, we assume that either $\kappa \leq 0$ or 
\[\frac{n^n}{4 \, (n-2)^{n-2}} \, m^2 \, \kappa^{n-2} < 1.\] 
Let us write $\{s > 0: 1 - m \, s^{2-n} - \kappa \, s^2 > 0\} = (\underline{s},\overline{s})$. The deSitter-Schwarzschild manifold is defined by $M = S^{n-1} \times (\underline{s},\overline{s})$ and 
\[g = \frac{1}{1 - m \, s^{2-n} - \kappa \, s^2} \, ds \otimes ds + s^2 \, g_{S^{n-1}}.\] 
A straightforward calculation shows that $(M,g)$ has constant scalar curvature $n(n-1)\kappa$. We note that the manifold $(M,g)$ can be isometrically embedded as a space-like slice in the $(n+1)$-dimensional deSitter-Schwarzschild space-time. In the special case $\kappa=0$, $(M,g)$ is the ordinary Schwarzschild manifold. 

\begin{corollary} 
\label{schwarzschild}
Suppose that $\Sigma$ is a closed, embedded, orientable hypersurface in the deSitter-Schwarzschild manifold with constant mean curvature. Then $\Sigma$ is a slice $S^{n-1} \times \{s\}$.
\end{corollary} 

In particular, we obtain a uniqueness theorem for surfaces of constant mean curvature in Schwarzschild space. We note, however, that the analogous result fails in the doubled Schwarzschild manifold. In fact, in a joint work with Michael Eichmair \cite{Brendle-Eichmair}, we have constructed small isoperimetric surfaces in the doubled Schwarzschild manifold which are located near a point on the horizon. These surfaces have constant mean curvature, but are not umbilic.

The assumptions of Theorem \ref{alexandrov} are also satisfied for the Reissner-Nordstrom spaces. The Reissner-Nordstrom manifold is defined by $M = S^{n-1} \times (\underline{s},\infty)$ and 
\[g = \frac{1}{1 - m \, s^{2-n} + q^2 \, s^{4-2n}} \, ds \otimes ds + s^2 \, g_{S^{n-1}}.\] 
Here, $m > 2q > 0$ are constants, and $\underline{s}$ is defined as the larger of the two solutions of the equation $1 - m \, s^{2-n} + q^2 \, s^{4-2n} = 0$.

\begin{corollary} 
\label{reissner.nordstrom}
Suppose that $\Sigma$ is a closed, embedded, orientable hypersurface in the Reissner-Nordstrom manifold with constant mean curvature. Then $\Sigma$ is a slice $S^{n-1} \times \{s\}$.
\end{corollary} 

We next state a variant of Theorem \ref{alexandrov}. To that end, we consider a function $h: [0,\bar{r}) \to \mathbb{R}$ which satisfies the following conditions:
\begin{itemize} 
\item[(H1')] $h(r) = r \, \varphi(r^2)$, where $\varphi: [0,\sqrt{\bar{r}}) \to \mathbb{R}$ is a smooth positive function satisfying $\varphi(0) = 1$.
\item[(H2')] $h'(r) > 0$ for all $r \in (0,\bar{r})$.
\item[(H3')] The function 
\[2 \, \frac{h''(r)}{h(r)} - (n-2) \, \frac{1 - h'(r)^2}{h(r)^2}\] 
is non-decreasing for $r \in (0,\bar{r})$.
\item[(H4')] We have 
\[\frac{h''(r)}{h(r)} + \frac{1-h'(r)^2}{h(r)^2} \neq 0\] 
for all $r \in (0,\bar{r})$.
\end{itemize}

\begin{theorem} 
\label{alexandrov.2}
Let $h: [0,\bar{r}) \to \mathbb{R}$ which satisfies the conditions (H1')--(H3'). Let us consider the ball $B_{\bar{r}}(0) \subset \mathbb{R}^n$ equipped with the Riemannian metric $g = dr \otimes dr + h(r)^2 \, g_{S^{n-1}}$. Moreover, let $\Sigma$ be a closed, embedded, orientable hypersurface in $(B_{\bar{r}}(0),g)$ with constant mean curvature. Then $\Sigma$ is umbilic. If, in addition, the condition (H4') holds, then $\Sigma$ is a geodesic sphere centered at the origin.
\end{theorem}

Note that the conditions (H1')--(H3') are satisfied for $h(r) = r$, for $h(r) = \sinh(r)$, and for $h(r) = \sin(r)$. Therefore, Theorem \ref{alexandrov.2} generalizes the Alexandrov theorems in Euclidean space, hyperbolic space, and the hemisphere.

As above, the condition (H3') is equivalent to saying that the scalar curvature of $g$ is non-increasing in $r$. On the other hand, if the scalar curvature is not a monotone function of $r$, it is possible to construct small spheres with constant mean curvature which are not umbilic:

\begin{theorem}[F.~Pacard, X.~Xu \cite{Pacard-Xu}]
\label{counterexample}
Consider a smooth metric $g$ on the ball $B_{\bar{r}}(0) \subset \mathbb{R}^n$ of the form $g = dr \otimes dr + h(r)^2 \, g_{S^{n-1}}$. Assume that the function 
\[2 \, \frac{h''(r)}{h(r)} - (n-2) \, \frac{1 - h'(r)^2}{h(r)^2}\] 
has a strict local extremum (either a local minimum or a local maximum) at some point $\hat{r} \in (0,\bar{r})$. Moroever, suppose that at each point on $\partial B_{\hat{r}}(0)$, the Ricci tensor of $g$ has two distinct eigenvalues. Then there exist small spheres with constant mean curvature which are not umbilic.
\end{theorem}

Theorem \ref{counterexample} is a direct consequence of Theorem 1.1 in \cite{Pacard-Xu}. To see this, let $\phi(\cdot,\lambda)$ be the function defined in \cite{Pacard-Xu}. Clearly, $\phi(\cdot,\lambda)$ is rotationally symmetric. By assumption, the scalar curvature (viewed as a function of $r$) attains a strict local extremum at the point $\hat{r} \in (0,\bar{r})$. Hence, if $\lambda$ is sufficiently small, we can find a number $r$ close to $\hat{r}$ with the property that every point on the coordinate sphere $\partial B_r(0)$ is a critical point of the function $\phi(\cdot,\lambda)$. Hence, given any point $p \in \partial B_r(0)$, we can find an $(n-1)$-sphere $\Sigma$ which has constant mean curvature $\frac{n-1}{\lambda}$ and is contained in a geodesic ball around $p$ of radius $\lambda \, (1+o(1))$. This surface $\Sigma$ cannot be umbilic: indeed, if $\Sigma$ were umbilic, then the Codazzi equations would imply that the normal vector to $\Sigma$ is an eigenvector of the Ricci tensor. Consequently, $\Sigma$ would be a geodesic sphere centered at the origin, which is impossible.

The proof of Theorem \ref{alexandrov} occupies Sections \ref{basic.properties} -- \ref{proof}. In Section \ref{basic.properties}, we define a potential function $f$ and a conformal vector field $X$, and study their basic properties. In particular, we derive an integral identity, which generalizes the classical Minkowski formula in Euclidean space. Moreover, using condition (H3), we show that $(\Delta f) \, g - D^2 f + f \, \text{\rm Ric} \geq 0$. This inequality plays a key role in the proof of Theorem \ref{alexandrov}. We note that Riemannian metrics with the property that $(\Delta f) \, g - D^2 f + f \, \text{\rm Ric} = 0$ are called static, and have been studied in connection with questions in general relativity (see e.g. \cite{Bartnik}, \cite{Corvino}, \cite{Reiris}). 

In Section \ref{inequality}, we prove a sharp inequality for hypersurfaces of positive mean curvature. This inequality is inspired by a classical inequality due to Heintze and Karcher \cite{Heintze-Karcher}. To state this inequality, we consider two cases. If $\Sigma$ is the boundary of a domain $\Omega$, we show that 
\begin{equation} 
\label{inequality.a}
(n-1) \int_\Sigma \frac{f}{H} \, d\mu \geq n \int_\Omega f \, d\text{\rm vol}. 
\end{equation} 
On the other hand, if we can find a domain $\Omega$ such that $\partial \Omega = \Sigma \cup (N \times \{0\})$, then we obtain 
\begin{equation} 
\label{inequality.b}
(n-1) \int_\Sigma \frac{f}{H} \, d\mu \geq n \int_\Omega f \, d\text{\rm vol} + h(0)^n \, \text{\rm vol}(N,g_N). 
\end{equation} 
Moreover, if equality holds in (\ref{inequality.a}) or (\ref{inequality.b}), then $\Sigma$ must be umbilic. In order to prove these inequalities, we consider the Riemannian metric $\hat{g} = \frac{1}{f^2} \, g$. This metric is conformal to the given metric $g$ and has an asymptotically hyperbolic end at $N \times \{0\}$. We then consider the level sets of the distance function $u(p) = d_{\hat{g}}(p,\Sigma)$, where the distance is computed using the conformal metric $\hat{g}$. We now study the quantity 
\[Q(t) = (n-1) \int_{\Sigma_t^*} \frac{f}{H} \, d\mu,\] 
where $\Sigma_t^*$ denotes the smooth part of the level set $\{u=t\}$, and $H$ denotes the mean curvature of $\Sigma_t^*$. Using the standard formula for the evolution of the mean curvature, we obtain a monotonicity formula for $Q(t)$. If $\Sigma$ is null-homologous, the inequality (\ref{inequality.a}) is a direct consequence of our monotonicity formula. On the other hand, if $\Sigma$ is homologous to the boundary $N \times \{0\}$, we show that 
\[\liminf_{t \to \infty} Q(t) \geq h(0)^n \, \text{\rm vol}(N,g_N).\] 
If we combine this inequality with the monotonicity formula for $Q(t)$, the inequality (\ref{inequality.b}) follows. The analysis of the limit of $Q(t)$ as $t \to \infty$ is a very subtle issue, as the level sets $\{u=t\}$ are not smooth in general. To overcome this obstacle, we use the approximation technique of Greene and Wu \cite{Greene-Wu1}, \cite{Greene-Wu2}.

In Section \ref{proof}, we combine the results from Sections \ref{basic.properties} and \ref{inequality} to conclude that $\Sigma$ is umbilic. This completes the proof of Theorem \ref{alexandrov}.

In Section \ref{subcase}, we explain how Corollary \ref{schwarzschild} and Corollary \ref{reissner.nordstrom} follow from Theorem \ref{alexandrov}.

Finally, in Section \ref{variant} we sketch the proof of Theorem \ref{alexandrov.2}. The proof of this result is very similar to the proof of Theorem \ref{alexandrov}, and we will indicate the necessary adaptations.

The author is very grateful to Professors Michael Eichmair, Jos\'e Espinar, and Brian White for discussions.

\section{Basic properties of $(M,g)$}

\label{basic.properties}

Let $(M,g)$ be a warped product manifold which satisfies the conditions (H1)--(H3). We define a smooth function $f: M \to \mathbb{R}$ and a vector field $X$ on $M$ by 
\[f = h'(r)\] 
and 
\[X = h(r) \, \frac{\partial}{\partial r}.\] 
It follows from (H2) that $f$ is a positive function on $N \times (0,\bar{r})$. Moreover, the condition (H1) implies that the function $f$ vanishes along $N \times \{0\}$, but the gradient of $f$ is non-zero along $N \times \{0\}$. \\

\begin{proposition}
\label{static}
The function $f$ satisfies the inequality 
\[(\Delta f) \, g - D^2 f + f \, \text{\rm Ric} \geq 0.\] 
\end{proposition}

\textbf{Proof.} 
Let $\{e_1,\hdots,e_{n-1}\}$ be a local orthonormal frame on $N$, so that $g_N(e_i,e_j) = \delta_{ij}$. It follows from Proposition 9.106 in \cite{Besse} that 
\begin{align*} 
&\text{\rm Ric}(e_i,e_j) = \text{\rm Ric}_N(e_i,e_j) - (h(r) \, h''(r) + (n-2) \, h'(r)^2) \, \delta_{ij}, \\ 
&\text{\rm Ric} \Big ( e_i,\frac{\partial}{\partial r} \Big ) = 0, \\ 
&\text{\rm Ric} \Big ( \frac{\partial}{\partial r},\frac{\partial}{\partial r} \Big ) = -(n-1) \, \frac{h''(r)}{h(r)}. 
\end{align*}
In other words, we have 
\begin{align*} 
\text{\rm Ric} 
&= \text{\rm Ric}_N - (h(r) \, h''(r) + (n-2) \, h'(r)^2) \, g_N \\ 
&- (n-1) \, \frac{h''(r)}{h(r)} \, dr \otimes dr. 
\end{align*}
On the other hand, the Hessian of $f$ is given by 
\[D^2 f = h(r) \, h'(r) \, h''(r) \, g_N + h'''(r) \, dr \otimes dr.\] 
This implies 
\begin{align*} 
(\Delta f) \, g - D^2 f 
&= (h(r)^2 \, h'''(r) + (n-2) \, h(r) \, h'(r) \, h''(r)) \, g_N \\ 
&+ (n-1) \, \frac{h'(r) \, h''(r)}{h(r)} \, dr \otimes dr. 
\end{align*}
Putting these facts together, we obtain 
\begin{align*} 
&(\Delta f) \, g - D^2 f + f \, \text{\rm Ric} \\ 
&= h'(r) \, (\text{\rm Ric}_N - (n-2)\rho \, g_N) \\ 
&+ \big ( h(r)^2 \, h'''(r) + (n-3) \, h(r) \, h'(r) \, h''(r) + (n-2) \, h'(r) \, (\rho - h'(r)^2) \big ) \, g_N. 
\end{align*} 
By assumption, we have 
\[\text{\rm Ric}_N - (n-2)\rho \, g_N \geq 0\] 
and 
\begin{align*} 
&h(r)^2 \, h'''(r) + (n-3) \, h(r) \, h'(r) \, h''(r) + (n-2) \, h'(r) \, (\rho - h'(r)^2) \\ 
&= \frac{1}{2} \, h(r)^3 \, \frac{d}{dr} \Big ( 2 \, \frac{h''(r)}{h(r)} - (n-2) \, \frac{\rho - h'(r)^2}{h(r)^2} \Big ) \geq 0. 
\end{align*} 
Putting these facts together, the assertion follows. \\

\begin{lemma}
\label{conformal}
The vector field $X$ satisfies $D_i X_j = f \, g_{ij}$.
\end{lemma}

\textbf{Proof.} 
Note that 
\[\mathscr{L}_X(dr) = d(X(r)) = d(h(r)) = h'(r) \, dr.\] 
This implies 
\begin{align*} 
\mathscr{L}_X g 
&= \mathscr{L}_X(dr) \otimes dr + dr \otimes \mathscr{L}_X(dr) + X(h(r)^2) \, g_N \\ 
&= 2h'(r) \, dr \otimes dr + 2h(r)^2 \, h'(r) \, g_N \\ 
&= 2h'(r) \, g. 
\end{align*}
Therefore, $\mathscr{L}_X g = 2f \, g$. Since $X$ is a gradient vector field, the assertion follows. \\

We next prove an analogue of the classical Minkowski formula (cf. \cite{Montiel}). Recall that $h''(0) > 0$ by condition (H1). By continuity, we can find a real number $r_1 \in (0,\bar{r})$ so that $h''(r) > 0$ for all $r \in [0,r_1]$.

\begin{proposition}
\label{minkowski}
Let $\Sigma$ be a closed orientable hypersurface in $(M,g)$. Then 
\begin{equation} 
\label{id.1}
\int_\Sigma H \, \langle X,\nu \rangle \, d\mu = (n-1) \int_\Sigma f \, d\mu. 
\end{equation}
Moreover, if $\Sigma$ is contained in the region $N \times (0,r_1)$, then we have the inequality 
\begin{equation} 
\label{id.2}
\int_\Sigma \frac{H}{f} \, \langle X,\nu \rangle \, d\mu \leq (n-1) \, \mu(\Sigma). 
\end{equation}
\end{proposition} 

\textbf{Proof.} 
Let us write $X = \nabla \psi$ for some real-valued function $\psi$. By Lemma \ref{conformal}, the Hessian of $\psi$ is given by $D^2 \psi = f \, g$. Hence, the Laplacian of the function $\psi|_\Sigma$ is given by 
\[\Delta_\Sigma \psi = \sum_{k=1}^{n-1} (D^2 \psi)(e_k,e_k) - H \, \langle \nabla \psi,\nu \rangle = (n-1)f - H \, \langle X,\nu \rangle,\] 
where $\{e_1,\hdots,e_{n-1}\}$ is an orthonormal basis of the tangent space to $\Sigma$. 
Therefore, we have 
\[(n-1) \int_\Sigma f \, d\mu - \int_\Sigma H \, \langle X,\nu \rangle \, d\mu = \int_\Sigma \Delta_\Sigma \psi \, d\mu = 0\] 
by the divergence theorem. This proves (\ref{id.1}). In order to prove the inequality (\ref{id.2}), we observe that 
\[(n-1) \, \mu(\Sigma) - \int_\Sigma \frac{H}{f} \, \langle X,\nu \rangle \, d\mu = \int_\Sigma \frac{1}{f} \, \Delta_\Sigma \psi \, d\mu = \int_\Sigma \frac{1}{f^2} \, \langle \nabla^\Sigma f,\nabla^\Sigma \psi \rangle \, d\mu.\] 
At each point in $N \times (0,r_1)$, the vector $\nabla f$ is a positive multiple of $\nabla \psi$. Hence, if $\Sigma$ is contained in the region $N \times (0,r_1)$, then we have $\langle \nabla^\Sigma f,\nabla^\Sigma \psi \rangle \geq 0$ at each point on $\Sigma$. From this, the assertion follows. \\

\section{A geometric inequality for mean-convex hypersurfaces}

\label{inequality}

We now consider a closed, embedded, orientable hypersurface $\Sigma$ in $(M,g)$. It is easy to see that the intersection number of $\Sigma$ with any closed loop is zero. Since $\Sigma$ is connected, the complement $M \setminus \Sigma$ has exactly two connected components. In particular, there is a unique connected component $\Omega$ of $M \setminus \Sigma$ with the property that $\Omega \subset N \times (0,\bar{r}-\delta)$ for some $\delta>0$. We either have $\partial \Omega = \Sigma$ or $\partial \Omega = \Sigma \cup (N \times \{0\})$. Let $\nu$ denote the outward-pointing unit normal to $\Sigma$. We will assume throughout this section that $\Sigma$ has positive mean curvature with respect to this choice of unit normal. 

It will be convenient to consider the conformally modified metric $\hat{g} = \frac{1}{f^2} \, g$. The manifold $(M,\hat{g})$ has an asymptotically hyperbolic end, which corresponds to the boundary $N \times \{0\}$. 

For each point $p \in \bar{\Omega}$, we denote by $u(p) = d_{\hat{g}}(p,\Sigma)$ the distance of $p$ from $\Sigma$ with respect to the metric $\hat{g}$. Moreover, we denote by $\Phi: \Sigma \times [0,\infty) \to \bar{\Omega}$ the normal exponential map with respect to $\hat{g}$. More precisely, for each point $x \in \Sigma$, the curve $t \mapsto \Phi(x,t)$ is a geodesic with respect to $\hat{g}$, and we have 
\[\Phi(x,0) = x, \qquad \frac{\partial}{\partial t} \Phi(x,t) \Big |_{t=0} = -f(x) \, \nu(x).\] 
Note that the geodesic $t \mapsto \Phi(x,t)$ has unit speed with respect to $\hat{g}$.

We next define 
\[A = \{(x,t) \in \Sigma \times [0,\infty): u(\Phi(x,t)) = t\}\] 
and 
\[A^* = \{(x,t) \in \Sigma \times [0,\infty): \text{\rm $(x,t+\delta) \in A$ for some $\delta>0$}\}.\] 
The definition of $A$ is analogous to the definition of the segment domain of a Riemannian manifold. Our next result follows from standard arguments (see e.g. \cite{Petersen}, pp.~139--141):

\begin{proposition}
The sets $A$ and $A^*$ have the following properties: \\ 
(i) If $(x,t_0) \in A$, then $(x,t) \in A$ for all $t \in [0,t_0]$. \\ 
(ii) The set $A$ is closed, and we have $\Phi(A) = \bar{\Omega}$. \\
(iii) The set $A^*$ is an open subset of $\Sigma \times [0,\infty)$, and the restriction $\Phi|_{A^*}$ is a diffeomorphism. 
\end{proposition}

For each $t \in [0,\infty)$, we define 
\[\Sigma_t^* = \Phi(A^* \cap (\Sigma \times \{t\})).\] 
Note that $\Sigma_t^*$ is a smooth hypersurface which is contained in the level set $\{u=t\}$. To fix notation, we denote by $H$ and $I\!I$ the mean curvature and second fundamental form of $\Sigma_t^*$ with respect to the metric $g$. 

\begin{proposition}
\label{mean.curvature}
The mean curvature of $\Sigma_t^*$ is positive and satisfies the differential inequality 
\[\frac{\partial}{\partial t} \Big ( \frac{f}{H} \Big ) \leq -\frac{1}{n-1} \, f^2.\] 
\end{proposition}

\textbf{Proof.} 
It is easy to see that $\frac{\partial}{\partial t} \Phi(x,t) = -f(\Phi(x,t)) \, \nu$, where $\nu = -\frac{\nabla u}{|\nabla u|}$ denotes the outward-pointing unit normal vector to $\Sigma_t^*$ with respect to the metric $g$. Hence, the mean curvature of $\Sigma_t^*$ satisfies the evolution equation 
\[\frac{\partial}{\partial t} H = \Delta_{\Sigma_t^*} f + (\text{\rm Ric}(\nu,\nu) + |I\!I|^2) \, f\] 
(cf. \cite{Huisken-Ilmanen}, equation (1.2)). Using Proposition \ref{static}, we obtain 
\begin{align*} 
\Delta_{\Sigma_t^*} f 
&= \sum_{k=1}^{n-1} (D^2 f)(e_k,e_k) - H \, \langle \nabla f,\nu \rangle \\ 
&= \Delta f - (D^2 f)(\nu,\nu) - H \, \langle \nabla f,\nu \rangle \\ 
&\geq -\text{\rm Ric}(\nu,\nu) \, f - H \, \langle \nabla f,\nu \rangle, 
\end{align*} 
where $\{e_1,\hdots,e_{n-1}\}$ is an orthonormal basis of the tangent space to $\Sigma_t^*$. Putting these facts together, we conclude that 
\[\frac{\partial}{\partial t} H \geq -H \, \langle \nabla f,\nu \rangle + |I\!I|^2 \, f.\] 
Moreover, we have 
\[\frac{\partial}{\partial t} f = -f \, \langle \nabla f,\nu \rangle.\] 
This implies 
\[\frac{\partial}{\partial t} \Big ( \frac{H}{f} \Big ) = \frac{1}{f} \, \frac{\partial}{\partial t} H - \frac{H}{f^2} \, \frac{\partial}{\partial t} f \geq |I\!I|^2 \geq \frac{1}{n-1} \, H^2\] 
at each point on $\Sigma_t^*$. Since the initial hypersurface $\Sigma$ has positive mean curvature, we conclude that the hypersurface $\Sigma_t^*$ has positive mean curvature for each $t \in [0,\infty)$. From this, the assertion follows.

\begin{corollary} 
\label{area.decreasing}
The function $t \mapsto \mu(\Sigma_t^*)$ is monotone decreasing.
\end{corollary}

\textbf{Proof.} 
Since $\Sigma_t^*$ has positive mean curvature, the area form on $\Sigma_t^*$ is monotone decreasing in $t$. Moreover, the sets $\{x \in \Sigma: (x,t) \in A^*\}$ become smaller as $t$ increases. From this, the assertion follows. \\

We next consider the quantity 
\[Q(t) = (n-1) \int_{\Sigma_t^*} \frac{f}{H} \, d\mu.\] 
It follows from Proposition \ref{mean.curvature} that the function $t \mapsto Q(t)$ is non-increasing. Moreover, we have the following estimate:\footnote{We note that Proposition \ref{monotonicity} can be extended to a more general setting. In fact, Michael Eichmair has pointed out that Proposition \ref{monotonicity} holds for any ambient manifold $(M,g)$ which satisfies the inequality $(\Delta f) \, g - D^2 f + f \, \text{\rm Ric} \geq 0$.}

\begin{proposition}
\label{monotonicity}
We have 
\[Q(0) - Q(\tau) \geq n \int_{\{u \leq \tau\}} f \, d\text{\rm vol}\] 
for all $\tau \in [0,\infty)$. 
\end{proposition}

\textbf{Proof.} 
Using Proposition \ref{mean.curvature}, we obtain 
\begin{align*} 
&\limsup_{h \searrow 0} \frac{1}{h} \, (Q(t) - Q(t-h)) \\ 
&\leq (n-1) \int_{\Sigma_t^*} \frac{\partial}{\partial t} \Big ( \frac{f}{H} \Big ) \, d\mu - (n-1) \int_{\Sigma_t^*} \frac{f}{H} \cdot fH \, d\mu \\ 
&\leq -n \int_{\Sigma_t^*} f^2 \, d\mu. 
\end{align*} 
Thus, we conclude that 
\begin{align*} 
Q(0) - Q(\tau) 
&\geq n \int_0^\tau \bigg ( \int_{\Sigma_t^*} f^2 \, d\mu \bigg ) \, dt \\ 
&= n \int_{\Phi(A^* \cap (\Sigma \times [0,\tau]))} f \, d\text{\rm vol} \\ 
&= n \int_{\{u \leq \tau\}} f \, d\text{\rm vol} 
\end{align*}
for all $\tau \in [0,\infty)$. \\

\begin{theorem} 
\label{heintze.karcher.1}
Assume that $\Sigma$ is null-homologous, so that $\partial \Omega = \Sigma$. Moreover, suppose that $\Sigma$ has positive mean curvature. Then 
\[(n-1) \int_\Sigma \frac{f}{H} \, d\mu \geq n \int_\Omega f \, d\text{\rm vol}.\] 
Moreover, if equality holds, then $\Sigma$ is umbilic.
\end{theorem}

\textbf{Proof.} 
By Proposition \ref{monotonicity}, we have 
\[Q(0) \geq n \int_{\{u \leq \tau\}} f \, d\text{\rm vol}\] 
for all $\tau \in [0,\infty)$. Passing to the limit as $\tau \to \infty$, we obtain 
\[Q(0) \geq n \int_\Omega f \, d\text{\rm vol},\] 
as claimed. \\

In the remainder of this section, we consider the case that $\Sigma$ is homologous to the boundary $N \times \{0\}$, so that $\partial \Omega = \Sigma \cup (N \times \{0\})$. Our goal is to analyze the asymptotics of $Q(\tau)$ when $\tau$ is very large. The key result is Proposition \ref{technical.ingredient}. The proof of this result is quite subtle, and relies on several lemmata: 

\begin{lemma}
\label{transversality.1}
Given any real number $\lambda \in (0,1)$, there exists a number $\tau_0>0$ with the following property: if $p$ is a point in $\{u \geq \tau_0\}$ and $\alpha$ is a unit-speed geodesic with respect to $\hat{g}$ such that $\alpha(0) = p$ and $\alpha(u(p)) \in \Sigma$, then $|\alpha'(0)| = f(p)$ and 
\[\Big \langle \frac{\partial}{\partial r},\alpha'(0) \Big \rangle \geq \lambda \, f(p).\] 
\end{lemma}

\textbf{Proof.} 
For abbreviation, let $c := h''(0) > 0$, so that $|df| = c$ along $N \times \{0\}$. By continuity, we can find a small number $r_0 \in (0,r_1)$ such that $\Sigma \subset N \times (r_0,\bar{r})$ and 
\[-f \, D^2 f + |df|^2 \, g \geq \lambda c^2 \, g\] 
on the set $N \times (0,r_0]$. Hence, the Hessian of the function $\frac{1}{f}$ with respect to $\hat{g}$ satisfies 
\begin{align*} 
\hat{D}^2 \Big ( \frac{1}{f} \Big ) 
&= D^2 \Big ( \frac{1}{f} \Big ) - \frac{2}{f^3} \, df \otimes df + \frac{1}{f^3} \, |df|^2 \, g \\ 
&= -\frac{1}{f^2} \, D^2 f + \frac{1}{f^3} \, |df|^2 \, g \\ 
&\geq \lambda c^2 \, \frac{1}{f} \, \hat{g} 
\end{align*} 
on the set $N \times (0,r_0]$. 

Let us choose $\tau_0$ sufficiently large so that $\{u \geq \tau_0-1\} \subset N \times (0,r_0]$ and 
\[c \, \Big ( 1 - \frac{f(p)}{h'(r_0) \, \sinh(\sqrt{\lambda} \, c)} \Big ) \geq \sqrt{\lambda} \, |\nabla f(p)|\] 
for all points $p \in \{u \geq \tau_0\}$. We claim that $\tau_0$ has the desired property. To verify this, we consider a point $p \in \{u \geq \tau_0\}$ and a unit-speed geodesic $\alpha$ with respect to $\hat{g}$ such that $\alpha(0) = p$ and $\alpha(u(p)) \in \Sigma$. Clearly, $\alpha(t) \in \{u \geq \tau_0-1\}$ for all $t \in [0,1]$. This implies $\alpha(t) \in N \times (0,r_0]$ for all $t \in [0,1]$. We now define $t_0 = \inf \{t \in [0,u(p)]: \alpha(t) \notin N \times (0,r_0]\}$. Clearly, $t_0 \geq 1$. Moreover, we have  
\[\frac{d^2}{dt^2} \Big ( \frac{1}{f(\alpha(t))} \Big ) \geq \lambda c^2 \, \frac{1}{f(\alpha(t))}\] 
for all $t \in [0,t_0]$. Integrating this differential inequality, we obtain 
\[\frac{1}{f(\alpha(t))} \geq \frac{1}{f(p)} \, \cosh(\sqrt{\lambda} \, c \, t) - \frac{1}{\sqrt{\lambda} \, c \, f(p)^2} \, \langle \nabla f(p),\alpha'(0) \rangle \, \sinh(\sqrt{\lambda} \, c \, t)\] 
for all $t \in [0,t_0]$. Putting $t=t_0$ and rearraning terms gives 
\[\langle \nabla f(p),\alpha'(0) \rangle \geq \sqrt{\lambda} \, c \, f(p) \, \Big ( \frac{\cosh(\sqrt{\lambda} \, c \, t_0)}{\sinh(\sqrt{\lambda} \, c \, t_0)} - \frac{f(p)}{h'(r_0) \, \sinh(\sqrt{\lambda} \, c \, t_0)} \Big ).\] 
Here, we have used the fact that $\alpha(t_0) \in N \times \{r_0\}$ and $f(\alpha(t_0)) = h'(r_0)$. On the other hand, we have 
\begin{align*} 
&c \, \Big ( \frac{\cosh(\sqrt{\lambda} \, c \, t_0)}{\sinh(\sqrt{\lambda} \, c \, t_0)} - \frac{f(p)}{h'(r_0) \, \sinh(\sqrt{\lambda} \, c \, t_0)} \Big ) \\ 
&\geq c \, \Big ( 1 - \frac{f(p)}{h'(r_0) \, \sinh(\sqrt{\lambda} \, c)} \Big ) \geq \sqrt{\lambda} \, |\nabla f(p)| 
\end{align*}
by our choice of $\tau_0$. Putting these facts together, we obtain 
\[\langle \nabla f(p),\alpha'(0) \rangle \geq \lambda \, f(p) \, |\nabla f(p)|.\] 
Since $\nabla f(p) = |\nabla f(p)| \, \frac{\partial}{\partial r}$, we conclude that 
\[\Big \langle \frac{\partial}{\partial r},\alpha'(0) \Big \rangle \geq \lambda \, f(p).\] 
This completes the proof of Lemma \ref{transversality.1}. \\

In the following, we fix a real number $\tau_1>0$ so that the conclusion of Lemma \ref{transversality.1} holds for $\lambda=\frac{1}{2}$. \\

\begin{lemma}
\label{transversality.2}
Suppose that $\gamma: [a,b] \to \{u \geq \tau_1\}$ is a smooth path satisfying $\big | \gamma'(s) + f(\gamma(s)) \, \frac{\partial}{\partial r} \big |_{\hat{g}} \leq \frac{1}{4}$ for all $s \in [a,b]$. Then 
\[u(\gamma(b)) - u(\gamma(a)) \geq \frac{1}{4} \, (b-a).\] 
\end{lemma}

\textbf{Proof.} 
It suffices to show that 
\begin{equation} 
\label{derivative}
\liminf_{h \searrow 0} \frac{1}{h} \, \big ( u(\gamma(s)) - u(\gamma(s-h)) \big ) \geq \frac{1}{4} 
\end{equation} 
for all $s \in (a,b]$. In order to verify (\ref{derivative}), we fix a real number $s_0 \in (a,b]$. For abbreviation, let $p = \gamma(s_0)$. Moreover, let $\alpha$ be a unit-speed geodesic with respect to $\hat{g}$ such that $\alpha(0) = p$ and $\alpha'(u(p)) \in \Sigma$. Then 
\[f(p) \, \Big \langle \frac{\partial}{\partial r},\alpha'(0) \Big \rangle_{\hat{g}} = \frac{1}{f(p)} \, \Big \langle \frac{\partial}{\partial r},\alpha'(0) \Big \rangle_g \geq \frac{1}{2}\] 
by our choice of $\tau_1$. Using the Cauchy-Schwarz inequality, we obtain 
\[-\langle \gamma'(s_0),\alpha'(0) \rangle_{\hat{g}} \geq f(p) \, \Big \langle \frac{\partial}{\partial r},\alpha'(0) \Big \rangle_{\hat{g}} - \Big | \gamma'(s_0) + f(p) \, \frac{\partial}{\partial r} \Big |_{\hat{g}} \geq \frac{1}{4}.\] 
Hence, it follows from the formula for the first variation of arclength that 
\[\liminf_{h \searrow 0} \frac{1}{h} \, \big ( u(\gamma(s_0)) - u(\gamma(s_0-h)) \big ) \geq -\langle \gamma'(s_0),\alpha'(0) \rangle_{\hat{g}} \geq \frac{1}{4}.\] 
This proves (\ref{derivative}), thereby completing the proof of Lemma \ref{transversality.2}. \\

In the next step, we approximate the function $u$ by smooth functions.

\begin{lemma}
\label{smoothing}
Given any real number $\tau \geq \tau_1+2$, there exists a sequence of smooth functions $u_j: \{\tau-1 < u < \tau+1\} \to \mathbb{R}$ with the following properties: 
\begin{itemize}
\item[(i)] The functions $u_j$ converge smoothly to $u$ away from the cut locus. More precisely, $u_j \to u$ in $C_{loc}^\infty(W)$, where $W = \Phi(A^* \cap (\Sigma \times (\tau-1,\tau+1)))$.
\item[(ii)] For each point $p \in \{\tau-1 < u < \tau+1\}$, we have $|u_j(p) - u(p)| \leq \frac{1}{j^2}$.
\item[(iii)] For all points $p,q \in \{\tau-1 < u < \tau+1\}$, we have $|u_j(p) - u_j(q)| \leq (1+\frac{1}{j}) \, d_{\hat{g}}(p,q)$.
\item[(iv)] If $\gamma: [a,b] \to \{\tau-1 < u < \tau+1\}$ is an integral curve of the vector field $-f \, \frac{\partial}{\partial r}$, then $u_j(\gamma(b)) - u_j(\gamma(a)) \geq \frac{1}{4} \, (b-a)$.
\item[(v)] We have $\hat{D}^2 u_j \leq K(\tau) \, \hat{g}$ at each point $p \in \{\tau-1 < u < \tau+1\}$. Here, $K(\tau)$ is a positive constant which may depend on $\tau$, but not on $j$.
\end{itemize}
\end{lemma}

\textbf{Proof.}
We employ the Riemannian convolution method of Greene and Wu (see \cite{Greene-Wu2}, p.~57). More precisely, we define 
\[u_j(p) = \int_{(T_p M,\hat{g})} G(|\xi|^2) \, u(\exp_p(\varepsilon_j \, \xi)) \, d\xi.\] 
Here, $\exp_p$ denotes the exponential map with respect to the metric $\hat{g}$ and $G: [0,\infty) \to [0,\infty)$ is a smooth function with compact support satisfying $\int_{\mathbb{R}^n} G(|\xi|^2) \, d\xi = 1$. Moreover, $\varepsilon_j$ is a sequence of positive real numbers which are chosen sufficiently small.

We claim that the sequence $u_j$ has the required properties. Properties (i) and (ii) are obvious. Property (iii) follows from the fact that $|u(p) - u(q)| \leq d_{\hat{g}}(p,q)$ for all points $p,q \in \bar{\Omega}$. Similarly, property (iv) is a consequence of Lemma \ref{transversality.2}. 

It remains to prove (v). We can find a positive real number $K(\tau)$ such that $\hat{D}^2 u \leq \frac{1}{2} \, K(\tau) \, \hat{g}$ at each point $p \in \{\tau-2 < u < \tau+2\}$, where the inequality is understood in the barrier sense. Results of Greene and Wu then imply that $\hat{D}^2 u_j \leq K(\tau) \, \hat{g}$ at each point $p \in \{\tau-1 < u < \tau+1\}$ (see \cite{Greene-Wu1}, p.~644, and \cite{Greene-Wu2}, p.~60). This completes the proof of Lemma \ref{smoothing}. \\

\begin{proposition} 
\label{technical.ingredient}
For $\tau \geq \tau_1+2$ we have 
\begin{equation} 
\label{technical.1}
\mu(\Sigma_\tau^*) \geq h(0)^{n-1} \, \text{\rm vol}(N,g_N) 
\end{equation}
and 
\begin{equation} 
\label{technical.2}
\int_{\Sigma_\tau^*} \frac{H}{f} \, \langle X,\nu \rangle \, d\mu \leq (n-1) \, \mu(\Sigma_\tau^*). 
\end{equation}
\end{proposition} 

\textbf{Proof.} 
Let us fix a real number $\tau \geq \tau_1+2$, and let $u_j$ be a sequence of smooth functions satisfying properties (i)--(v) in Lemma \ref{smoothing}. The statement (iii) implies that 
\[f \, |du_j|_g = |du_j|_{\hat{g}} \leq 1+\frac{1}{j}.\] 
Moreover, it follows from (iv) that 
\[-f \, \frac{\partial}{\partial r} u_j \geq \frac{1}{4}.\] 
In particular, we have $f \, |du_j|_g \geq \frac{1}{4}$ and $\langle X,-\nabla u_j \rangle > 0$.

Using the co-area formula, we obtain 
\begin{align*} 
\int_{\tau+\frac{1}{j^2}}^{\tau+\frac{1}{j}-\frac{1}{j^2}} \mu(\{u_j = t\}) \, dt 
&= \int_{\{\tau+\frac{1}{j^2} \leq u_j \leq \tau+\frac{1}{j}-\frac{1}{j^2}\}} |\nabla u_j| \, d\text{\rm vol} \\ 
&\leq \Big ( 1+\frac{1}{j} \Big ) \int_{\{\tau \leq u \leq \tau+\frac{1}{j}\}} \frac{1}{f} \, d\text{\rm vol}, 
\end{align*} 
where the volume form is taken with respect to the metric $g$. Moreover, using Corollary \ref{area.decreasing}, we obtain 
\begin{align*} 
\int_{\{\tau \leq u \leq \tau+\frac{1}{j}\}} \frac{1}{f} \, d\text{\rm vol} 
&= \int_{\Phi(A^* \cap (\Sigma \times [\tau,\tau+\frac{1}{j}]))} \frac{1}{f} \, d\text{\rm vol} \\ 
&= \int_\tau^{\tau+\frac{1}{j}} \mu(\Sigma_t^*) \, dt \\ 
&\leq \frac{1}{j} \, \mu(\Sigma_\tau^*). 
\end{align*}
Putting these facts together, we conclude that 
\[\int_{\tau+\frac{1}{j^2}}^{\tau+\frac{1}{j}-\frac{1}{j^2}} \mu(\{u_j = t\}) \, dt \leq \frac{1}{j} \, \Big ( 1 + \frac{1}{j} \Big ) \, \mu(\Sigma_\tau^*).\] 
Therefore, we can find a real number $t_j \in [\tau+\frac{1}{j^2},\tau+\frac{1}{j}-\frac{1}{j^2}]$ such that 
\begin{equation} 
\label{choice.of.tj}
\mu(\{u_j = t_j\}) \leq \frac{j+1}{j-2} \, \mu(\Sigma_\tau^*). 
\end{equation}
Let us denote the level set $\{u_j = t_j\}$ by $S_j$. Clearly, $S_j$ is a smooth (possibly disconnected) hypersurface without boundary. Using (\ref{choice.of.tj}), we obtain 
\begin{equation} 
\label{upper.bound.for.area}
\limsup_{j \to \infty} \mu(S_j) \leq \mu(\Sigma_\tau^*). 
\end{equation}
It follows from the intermediate value theorem that every integral curve of the vector field $\frac{\partial}{\partial r}$ intersects $S_j$ at least once. Therefore, we have $\mu(S_j) \geq h(0)^{n-1} \, \text{\rm vol}(N,g_N)$. Passing to the limit as $j \to \infty$, we conclude that $\mu(\Sigma_\tau^*) \geq h(0)^{n-1} \, \text{\rm vol}(N,g_N)$. This proves (\ref{technical.1}).

It remains to verify the inequality (\ref{technical.2}). The outward-pointing unit normal vector to the hypersurface $S_j$ is given by $-\frac{\nabla u_j}{|\nabla u_j|}$. Moreover, the mean curvature of $S_j$ is given by 
\[H_{S_j} = -\frac{1}{|\nabla u_j|} \, \Big ( \Delta u_j - \frac{(D^2 u_j)(\nabla u_j,\nabla u_j)}{|\nabla u_j|^2} \Big ).\] 
We will denote by $S_j^+$ the set of all points on $S_j$ where the mean curvature $H_{S_j}$ is positive.

In view of property (i) above, the surfaces $S_j$ converge to $\Sigma_\tau^*$ in $C_{loc}^\infty$ away from the cut locus. Since $\Sigma_\tau^*$ has positive mean curvature, we have 
\begin{equation} 
\label{lower.bound.for.area}
\liminf_{j \to \infty} \mu(S_j^+) \geq \mu(\Sigma_\tau^*) 
\end{equation}
and 
\begin{equation} 
\label{integral.1}
\liminf_{j \to \infty} \int_{S_j^+} \frac{H_{S_j}}{f} \, \Big \langle X,-\frac{\nabla u_j}{|\nabla u_j|} \Big \rangle \, d\mu \geq \int_{\Sigma_\tau^*} \frac{H}{f} \, \langle X,\nu \rangle \, d\mu. 
\end{equation}
Combining (\ref{upper.bound.for.area}) and (\ref{lower.bound.for.area}), we obtain 
\begin{equation} 
\label{bad.set.has.small.measure}
\limsup_{j \to \infty} \mu(S_j \setminus S_j^+) = 0. 
\end{equation}
On the other hand, it follows from property (v) above that $D^2 u_j \leq L(\tau) \, g$ for some positive constant $L(\tau)$. Since $f \, |\nabla u_j| \geq \frac{1}{4}$, we conclude that $H_{S_j} \geq -\Lambda(\tau)$ for some positive constant $\Lambda(\tau)$. Note that the constants $L(\tau)$ and $\Lambda(\tau)$ may depend on $\tau$, but not on $j$. Using (\ref{bad.set.has.small.measure}), we obtain 
\begin{align} 
\label{integral.2}
&\liminf_{j \to \infty} \int_{S_j \setminus S_j^+} \frac{H_{S_j}}{f} \, \Big \langle X,-\frac{\nabla u_j}{|\nabla u_j|} \Big \rangle \, d\mu \notag \\ 
&= \liminf_{j \to \infty} \int_{S_j \setminus S_j^+} \frac{H_{S_j}+\Lambda(\tau)}{f} \, \Big \langle X,-\frac{\nabla u_j}{|\nabla u_j|} \Big \rangle \, d\mu \geq 0. 
\end{align} 
Adding (\ref{integral.1}) and (\ref{integral.2}) gives 
\[\liminf_{j \to \infty} \int_{S_j} \frac{H_{S_j}}{f} \, \Big \langle X,-\frac{\nabla u_j}{|\nabla u_j|} \Big \rangle \, d\mu \geq \int_{\Sigma_\tau^*} \frac{H}{f} \, \langle X,\nu \rangle \, d\mu.\] 
On the other hand, applying Proposition \ref{minkowski} to the hypersurface $S_j$ yields 
\begin{align*} 
\limsup_{j \to \infty} \int_{S_j} \frac{H_{S_j}}{f} \, \Big \langle X,-\frac{\nabla u_j}{|\nabla u_j|} \Big \rangle \, d\mu 
&\leq (n-1) \limsup_{j \to \infty} \mu(S_j) \\ 
&\leq (n-1) \, \mu(\Sigma_\tau^*). 
\end{align*}
Thus, we conclude that 
\[\int_{\Sigma_\tau^*} \frac{H}{f} \, \langle X,\nu \rangle \, d\mu \leq (n-1) \, \mu(\Sigma_\tau^*).\] 
This completes the proof of Proposition \ref{technical.ingredient}. \\

\begin{corollary} 
\label{asymptotics.of.Q}
Let $\lambda \in (0,1)$ be given. Then we have 
\[(n-1) \int_{\Sigma_\tau^*} \frac{f}{H} \, d\mu \geq \lambda \, h(0)^n \, \text{\rm vol}(N,g_N)\] 
if $\tau$ is sufficiently large.
\end{corollary}

\textbf{Proof.} 
It follows from Lemma \ref{transversality.1} that $\inf_{\Sigma_\tau^*} \langle \frac{\partial}{\partial r},\nu \rangle \geq \lambda$ if $\tau$ is sufficiently large. This implies that $\inf_{\Sigma_\tau^*} \langle X,\nu \rangle \geq \lambda \, h(0)$ if $\tau$ is sufficiently large. Using Proposition \ref{technical.ingredient}, we obtain 
\[\lambda \, h(0) \int_{\Sigma_\tau^*} \frac{H}{f} \, d\mu \leq \int_{\Sigma_\tau^*} \frac{H}{f} \, \langle X,\nu \rangle \, d\mu \leq (n-1) \, \mu(\Sigma_\tau^*).\] 
This implies 
\begin{align*} 
(n-1) \int_{\Sigma_\tau^*} \frac{f}{H} \, d\mu 
&\geq (n-1) \, \mu(\Sigma_\tau^*)^2 \, \bigg ( \int_{\Sigma_\tau^*} \frac{H}{f} \, d\mu \bigg )^{-1} \\ 
&\geq \lambda \, h(0) \, \mu(\Sigma_\tau^*) \\ 
&\geq \lambda \, h(0)^n \, \text{\rm vol}(N,g_N). 
\end{align*}
This completes the proof of Corollary \ref{asymptotics.of.Q}. \\

\begin{theorem} 
\label{heintze.karcher.2}
Assume that $\Sigma$ is homologous to the boundary $N \times \{0\}$, so that $\partial \Omega = \Sigma \cup (N \times \{0\})$. Moreover, suppose that $\Sigma$ has positive mean curvature. Then 
\[(n-1) \int_\Sigma \frac{f}{H} \, d\mu \geq n \int_\Omega f \, d\text{\rm vol} + h(0)^n \, \text{\rm vol}(N,g_N).\] 
Moreover, if equality holds, then $\Sigma$ is umbilic.
\end{theorem}

\textbf{Proof.} 
By Proposition \ref{monotonicity}, we have  
\[Q(0) - Q(\tau) \geq \int_{\{u \leq \tau\}} f \, d\text{\rm vol}\] 
for all $\tau \in [0,\infty)$. Moreover, we have 
\[\liminf_{\tau \to \infty} Q(\tau) \geq h(0)^n \, \text{\rm vol}(N,g_N)\] 
by Corollary \ref{asymptotics.of.Q}. Putting these facts together, we conclude that 
\[Q(0) \geq n \int_\Omega f \, d\text{\rm vol} + h(0)^n \, \text{\rm vol}(N,g_N),\] 
as claimed. \\

\section{Proof of Theorem \ref{alexandrov}}

\label{proof}

In this section, we give the proof of Theorem \ref{alexandrov}. As above, we assume that $(M,g)$ is a warped product manifold satisfying conditions (H1)--(H3). Let $\Sigma$ be a closed, embedded, orientable hypersurface in $(M,g)$ with constant mean curvature. It follows from (H2) that the slice $N \times \{r\}$ has positive mean curvature for each $r \in (0,\bar{r})$. Using the maximum principle, we conclude that the mean curvature of $\Sigma$ is strictly positive. By Proposition \ref{minkowski}, we have 
\[(n-1) \int_\Sigma f \, d\mu = \int_\Sigma H \, \langle X,\nu \rangle \, d\mu.\] 
Since $H$ is constant, we obtain 
\[(n-1) \int_\Sigma \frac{f}{H} \, d\mu = \int_\Sigma \langle X,\nu \rangle \, d\mu.\] 
We now distinguish two cases: 

\textit{Case 1:} Suppose first that $\Sigma$ is null-homologous, so that $\partial \Omega = \Sigma$. Using Lemma \ref{conformal} and the divergence theorem, we obtain  
\begin{align*} 
(n-1) \int_\Sigma \frac{f}{H} \, d\mu 
&= \int_\Sigma \langle X,\nu \rangle \, d\mu \\ 
&= \int_\Omega \text{\rm div} \, X \, d\text{\rm vol} \\ 
&= n \int_\Omega f \, d\text{\rm vol}. 
\end{align*}
Therefore, it follows from Theorem \ref{heintze.karcher.1} that $\Sigma$ is umbilic. 

\textit{Case 2:} We now assume that $\Sigma$ is homologous to the boundary $N \times \{0\}$, so that $\partial \Omega = \Sigma \cup (N \times \{0\})$. In this case, we have 
\begin{align*} 
(n-1) \int_\Sigma \frac{f}{H} \, d\mu 
&= \int_\Sigma \langle X,\nu \rangle \, d\mu \\ 
&= \int_\Omega \text{\rm div} \, X \, d\text{\rm vol} + h(0)^n \, \text{\rm vol}(N,g_N) \\ 
&= n \int_\Omega f \, d\text{\rm vol} + h(0)^n \, \text{\rm vol}(N,g_N). 
\end{align*}
Thus, Theorem \ref{heintze.karcher.2} implies that $\Sigma$ is umbilic. \\

Finally, let us assume that the condition (H4) is satisfied. In this case, we claim that $\Sigma$ is a slice $N \times \{r\}$. We have already shown that the second fundamental form of $\Sigma$ is a constant multiple of the metric. Using the Codazzi equations, we deduce that $R(e_i,e_j,e_k,\nu) = 0$, where $\{e_1,\hdots,e_{n-1}\}$ is an orthonormal basis for the tangent space of $\Sigma$. In particular, $\text{\rm Ric}(e_j,\nu) = \sum_{i=1}^{n-1} R(e_i,e_j,e_i,\nu) = 0$. Therefore, $\nu$ must be an eigenvector of the Ricci tensor of $(M,g)$. On the other hand, the condition (H4) implies that the smallest eigenvalue of the Ricci tensor is equal to $-(n-1) \, \frac{h''(r)}{h(r)}$; moreover, the corresponding eigenspace is one-dimensional and is spanned by the vector $\frac{\partial}{\partial r}$. Hence, at each point on $\Sigma$, the unit normal vector $\nu$ is either parallel or orthogonal to the vector $\frac{\partial}{\partial r}$. However, there is at least one point on $\Sigma$ where $\nu$ is parallel to $\frac{\partial}{\partial r}$. Thus, $\nu$ is parallel to $\frac{\partial}{\partial r}$ at each point on $\Sigma$, and $\Sigma$ is a slice $N \times \{r\}$.

\section{Application to the deSitter-Schwarzschild and Reissner-Nordstrom manifolds}

\label{subcase}

In this section, we describe how Corollary \ref{schwarzschild} and Corollary \ref{reissner.nordstrom} follow from Theorem \ref{alexandrov}. Let us consider the product $M = N \times (\underline{s},\overline{s})$ equipped with a metric of the form 
\[g = \frac{1}{\omega(s)} \, ds \otimes ds + s^2 \, g_{S^{n-1}}.\] 
Here, $\overline{s} > \underline{s} > 0$, and $\omega$ is a smooth function defined on the interval $[\underline{s},\overline{s})$. 

To bring the metric into the form (\ref{metric}), we define a continuous function $F: [\underline{s},\overline{s}) \to \mathbb{R}$ by $F'(s) = \frac{1}{\sqrt{\omega(s)}}$ and $F(\underline{s}) = 0$. Using the subsitution $r = F(s)$, the metric can be rewritten as 
\[g = dr \otimes dr + h(r)^2 \, g_{S^{n-1}},\] 
where $h: [0,F(\overline{s})) \to [\underline{s},\overline{s})$ denotes the inverse of the function $F$. A straightforward calculation gives 
\[h'(r) = \sqrt{\omega(s)}\] 
and 
\[h''(r) = \frac{1}{2} \, \omega'(s)\] 
where $s = h(r)$. Hence, the conditions (H1)--(H4) are equivalent to the following set of conditions: 
\begin{itemize}
\item $\omega(\underline{s}) = 0$ and $\omega'(\underline{s}) > 0$.
\item The function 
\[\frac{\omega'(s)}{s} - (n-2) \, \frac{\rho - \omega(s)}{s^2}\] 
is non-decreasing for $s \in (\underline{s},\overline{s})$.
\item We have 
\[\frac{\omega'(s)}{2s} + \frac{\rho-\omega(s)}{s^2} > 0\] 
for all $s \in (\underline{s},\overline{s})$.
\end{itemize}
Note that $\omega(s) = 1 - m \, s^{2-n} - \kappa \, s^2$ for the deSitter-Schwarzschild manifold, and $\omega(s) = 1 - m \, s^{2-n} + q^2 \, s^{4-2n}$ for the Reissner-Nordstrom manifold. Moreover, we have $\rho=1$ in both cases. It is straightforward to verify that the conditions above are satisfied. Thus, we can apply Theorem \ref{alexandrov} to these manifolds.


\section{Proof of Theorem \ref{alexandrov.2}}

\label{variant}

In this final section, we sketch the proof of Theorem \ref{alexandrov.2}. Let $h: [0,\bar{r}) \to \mathbb{R}$ be a smooth function which satisfies the conditions (H1')--(H3'). We define a Riemannian metric $g$ on the ball $B_{\bar{r}}(0) \subset \mathbb{R}^n$ by $g = dr \otimes dr + h(r)^2 \, g_{S^{n-1}}$. The condition (H1') implies that $g$ is smooth. As above, we define $f = h'(r)$ and $X = h(r) \, \frac{\partial}{\partial r}$. Note that $f$ is a smooth positive function defined on $B_{\bar{r}}(0) \subset \mathbb{R}^n$, and $X$ is a smooth vector field.

We now assume that $\Sigma$ is a closed, embedded, orientable hypersurface in $(B_{\bar{r}}(0),g)$ with constant mean curvature. Moreover, let $\Omega \subset B_{\bar{r}}(0)$ denote the domain enclosed by $\Sigma$. By assumption, the coordinate spheres $\partial B_r(0)$ have positive mean curvature for each $r \in (0,\bar{r})$. This implies that the mean curvature of $\Sigma$ is strictly positive. 

For each point $p \in \bar{\Omega}$, we denote by $u(p)$ the distance of $p$ from $\Sigma$ with respect to the metric $\hat{g} = \frac{1}{f^2} \, g$. It follows from Proposition \ref{monotonicity} that 
\[(n-1) \int_\Sigma \frac{f}{H} \, d\mu \geq n \int_{\{u \leq \tau\}} f \, d\text{\rm vol}\] 
for each $\tau \in [0,\infty)$. Passing to the limit as $\tau \to \infty$, we obtain 
\begin{equation} 
\label{inequality.2}
(n-1) \int_\Sigma \frac{f}{H} \, d\mu \geq n \int_\Omega f \, d\text{\rm vol}. 
\end{equation}
Moreover, if equality holds, then $\Sigma$ is umbilic. 

On the other hand, it follows from Proposition \ref{minkowski} that 
\[(n-1) \int_\Sigma f \, d\mu = \int_\Sigma H \, \langle X,\nu \rangle \, d\mu.\] 
Since $H$ is constant, we conclude that 
\[(n-1) \int_\Sigma \frac{f}{H} \, d\mu = \int_\Sigma \langle X,\nu \rangle \, d\mu = \int_\Omega \text{\rm div} \, X \, d\text{\rm vol} = n \int_\Omega f \, d\text{\rm vol}.\] 
Therefore, equality holds in (\ref{inequality.2}). Thus, $\Sigma$ is umbilic. If the condition (H4') holds, then the Ricci tensor of $g$ has two distinct eigenvalues. Hence, we can argue as in Section \ref{proof} to conclude that $\Sigma$ is a geodesic sphere centered at the origin.

\end{document}